\documentclass[11pt,reqno]{amsart}
\usepackage[pdftex, colorlinks=true, citecolor=blue, linkcolor=blue,
linktocpage=true]{hyperref}
\usepackage[utf8]{inputenc}
\usepackage{amsmath,amssymb,amsthm}
\usepackage{listings}
\usepackage{hyperref}

\usepackage[a4paper,margin=1.2in]{geometry}
\usepackage{hyperref}
\usepackage{graphicx}
\usepackage{tikz}
\usepackage{tikz-cd}

\newtheorem{fact}{Fact}[section]
\newtheorem{theorem}[fact]{Theorem}
\newtheorem*{theorem*}{Theorem}

\newtheorem{corollary}[fact]{Corollary}
\newtheorem{remark}[fact]{Remark}

\newtheorem{example}[fact]{Example}

\def\Q{{\mathbb Q}}
\def\G{{\mathbb G}}
\def\P{{\mathbb P}}
\def\T{{\mathbb T}}
\def\C{{\mathbb C}}

\def\ccc{{\mathfrak C}}
\def\cc{{\mathfrak c}}
\def\Grass{{\rm Grass}}
\def\Hom{{\rm Hom}}
\def\mepsilon{{\bf r}}
\def\fgf{f}
\def\ff{\phi}

\title[Four paths to the elliptic genus]{\huge$\text{\normalfont Four paths from
birational geometry}$\\[0.5cm] $\text{\normalfont to the elliptic genus}$}
\author{Andrzej Weber}\address{Institute of Mathematics, University of Warsaw,
Poland}
\email{aweber@mimuw.edu.pl}
\thanks{A.W.~supported by Polish National Science Center grant number
2022/47/B/ST1/01896
\hfill\break\phantom{aa}
}

\begin{document}

\maketitle

\begin{abstract}
The article presents four reasons why the elliptic genus is the most general
characteristic class that admits a generalization to singular spaces.
We prove that the elliptic characteristic class (with an additional factor) is
essentially the only characteristic class invariant under certain modifications,
such as the Atiyah flop, Grassmannian flops, and modifications of Bott-Samelson
resolutions.
This result confirms and extends Totaro's result concerning the cobordism ring
modulo classical flops. However, our approach is based on local calculus in
equivariant cohomology.
\end{abstract}
\tableofcontents

\section{Introduction}
Our goal is to address the following question:

\begin{center}\emph{Which characteristic classes can be reasonably extended to
singular varieties?}\end{center}

\noindent Perhaps we wish to extend not to all types of singular varieties, but
at least to a broad and meaningful family.

For smooth varieties, characteristic classes are typically defined via the
tangent bundle. However, singular varieties lack a globally defined tangent
bundle, which requires alternative constructions. Since the seventies of the
last century, there have appeared various constructions of characteristic
classes for complex singular varieties: first, the Chern-Schwartz-MacPherson
classes, and later, the Todd-Baum-Fulton-MacPherson classes. Also there purely
topological construction of L-classes followed from Goresky-Mac\-Pher\-son
theory of intersection cohomology. In the 21st century motivic/Hirzebruch Chern
classes were constructed by Brasselet-Sch\"urmann-Yokura and finally the
Borisov-Libgober elliptic classes. See the survey \cite{SYsurvey}. The most
straightforward way to define a characteristic class for a singular variety is
to find a resolution of singularities, compute its characteristic class, and
apply the appropriate push-forward to the singular variety. In general, the
result depends on the resolution chosen, and the contributions of singular loci
must be accounted for as correction terms.

There are two common approaches to studying singularities. The first involves
decomposing singular spaces into manifolds, called strata, and summing
invariants that depend on each stratum and its associated normal data. In this
work, we follow a different route: we replace the singular variety with a smooth
one via resolution and then combine the invariants of the resolution with those
of the components of the exceptional divisor.
Any two resolutions of a singular variety $\psi_i:X_i\to Y$ are dominated by a
third one
$$\begin{matrix} &_{\varphi_1}&X_3&_{\varphi_2}\\
&\swarrow&&\searrow&\\
X_1&&&&X_2\\
&\searrow&&\swarrow\\
&^{\psi_1}&Y&^{\psi_2}\end{matrix}\,.$$
Therefore it is enough to understand invariants of smooth varieties and their
behavior with respect to birational morphisms. By weak factorization theorem,
\cite{Wlod}, the morphisms $\varphi_i$ can be assumed to be blow-ups in smooth
centers.

By a {\it characteristic class}, we mean an element of the cohomology ring
arising from the Hirzebruch formalism of multiplicative characteristic classes.
To start we are interested in those characteristic classes $\cc$ that satisfy
the following birational invariance property: for any birational proper morphism
of complex algebraic varieties $\varphi\colon X' \to X$, we have
\begin{equation}\label{naive}
\varphi_* \cc(X')=\cc(X)\,.
\end{equation}
Such classes are rare, as is well known and will also be demonstrated through a
straightforward computation in \S\ref{baby}. This scarcity suggests two possible
approaches: either one must take into account the invariants of the exceptional
locus, or restrict attention to a certain class of admissible resolutions. Both
approaches will be discussed.

Our main result can be summarized as follows: we focus on the following settings

\begin{enumerate}
\item In the first setting, instead of demanding the naive equality
\eqref{naive} we correct the formula by the contributions depending on the
exceptional divisor. In fact, instead of characteristic classes of singular
varieties we study relative characteristic classes of pairs consisting of a
variety with a divisor.
As in \cite{BL} we assume that the pair has log-canonical singularities (see
Remark \ref {logcan}).
We single out the condition imposed by the invariance with respect to the
blow-up in a center of codimension two.

\item In the second setting, we compare the characteristic classes of
Bott-Samelson resolutions of a given Schubert variety. Here, we also consider
the relative characteristic classes of pairs. We examine the consequences of the
braid relation that connect two different resolutions.

\item As it is common in birational geometry, we focus on a smaller class of
singular varieties: those admiting a {\it crepant} resolution. We test what are
the consequences of the Atiyah flop, i.e. the basic example of two,
nonequivalent resolutions of the affine cone over $\P^1\times \P^1$.

\item Finally, we consider symplectic singularities. The best-known examples are
the closures of nilpotent orbits in semisimple Lie algebras. Once again, we
focus on the simplest possible cases of such singular varieties and examine the
consequences of choosing two non-equivalent symplectic
resolutions.\end{enumerate}

\begin{theorem*}[Theorems \ref{pierwszy}, \ref{drugi} \ref{AtFlop},
\ref{czwarty}] In each of the cases:
\begin{enumerate}
\item pairs, with at worst log-canonical singularities

\item Schubert varieties with Bott-Samelson resolutions,

\item singularities admitting crepant resolutions

\item symplectic singularities.
\end{enumerate}
the only (up to rescaling) characteristic class that is invariant under
birational morphisms is the elliptic class multiplied by the factor $e^{\alpha
c_1(X)}$ for
some $\alpha \in \C$.\end{theorem*}

The true meaning of our theorem is that, by restricting the classes of singularities, we do not obtain any new invariant characteristic classes other than the elliptic class, up to a minor modification.
\medskip

To place our result in a broader context, let us first recall the significant
work of Totaro~\cite{Totaro}, who showed that the ideal in the cobordism ring
$MU \otimes \Q$ consisting of almost complex manifolds with vanishing elliptic
genus is generated by the difference $X_1-X_2$, where $X_1$ and $X_2$ are
related by a classical flop. This result is based on the previous study of the
elliptic genus for singular varieties by Borisov and Libgober~\cite{BL}.

\begin{theorem*}[{\cite[Theorem 4.1]{Totaro}}] Let $I$ be the ideal in the
complex bordism ring $MU_* \otimes \Q$
which is additively generated by differences $X_1 - X_2$, where $X_1$ and $X_2$
are
smooth projective varieties related by a classical flop. Then the complex
elliptic genus
induces an isomorphism $$\Q[x_1,x_2,x_3,x_4]\to(MU_* \otimes \Q)/I \,,$$ where
$x_i=[\P^i]$. \end{theorem*}

It remains to note that the algebra in question is
naturally isomorphic to the algebra of quasi-Jacobi forms of depth $(k, 0)$, $k
\geq 0$, {\cite[Theorem 2.12]{Lib}}.
Our result is similar in spirit, though it relies solely on a local approach and
is formulated in terms of characteristic classes. The proof demonstrates that
the constraints on admissible classes arise from comparing pairs of resolutions
of particular singularities. That is consistent with the theorem quoted above.
While the case of Atiyah flop can be deduced from Totaro's theorem, what appears
to be a new result, is the case of symplectic singularities. We show that
demanding invariance of a
characteristic class under two independent Grassmannian flops determines the
elliptic classes as well.

The idea of identifying conditions under which a general characteristic class
remains invariant under specific modifications appears in~\cite{BrEv}, where the
authors studied fundamental classes of Schubert varieties in generalized
cohomology theories. There, it was shown that essentially only homological and
$K$-theoretic fundamental classes are well-defined. Applying the generalized
Riemann-Roch theorem \cite[\S42]{FFG}, the above fact translates into a
statement about characteristic classes: only the Todd class (and its
deformations) admits a generalization to Schubert varieties.

The scarcity of characteristic classes which are functorial on the nose for
birational morphism can be overcome by modifying the notion of functoriality. By
incorporating nontrivial contributions from the singular fibers into the
definition of the pushforward, one can broaden the family of characteristic
classes available for singular varieties. Most notably, the
Chern-Schwartz-MacPherson classes should be mentioned. Furthermore, as shown
in~\cite{BSY}, Hirzebruch classes also admit generalizations to singular
varieties. These can be interpreted as fundamental classes in hyperbolic
$K$-theory~\cite{LenartZa}. A wide range of related characteristic classes for
Schubert varieties (more precisely, Schubert cells) is presented in a series of
papers by the authors of~\cite{AMSS}. Various application to representation
theory have been found.

Elliptic classes of varieties with at worst Kawamata log-terminal singularities
were introduced by Borisov and Libgober in~\cite{BL}. Moreover, defining
characteristic classes for pairs turned out to be essential. As part of their
construction, they defined a characteristic class for a smooth variety together
with a $\Q$-divisor. This characteristic class was associated with a power
series given by the Taylor expansion of the Jacobi theta function.
Our guiding question regarding the most general form of characteristic classes
for pairs was simply: Can the expanion of the theta function be replaced by
another power series?
The answer we obtained is that such a replacement is indeed possible:
an additional factor can be introduced.

Our proof is based on elementary calculations once we adopt the perspective of
equivariant cohomology. It is striking that simple and concrete
computations can lead to such general result.
The basic (and practically the only)
tool is the Atiyah-Bott, Berline-Vergne localization theorem (see Theorem
\ref{ABtheorem}).

As and addendum to Theorem \ref{pierwszy} we show that a function $\fgf(x)$ defines a characteristic class of pairs, which is invariant with respect to blowups, if and only if $u(x)=\log\left(\frac {f(x)}x\right)''$ satisfies the differential equation
\begin{equation}\label{mojerownanie}u''(x)= 12\left(\frac{u(x)}{x^2}+u(x)^2\right)+C\,.\end{equation}
Since the Borisov-Libgober elliptic genus is preserved by blow-ups, there is a solution of the form \begin{equation}\label{weierstarss}u(x)=\sum_{k=1}^\infty
\frac{\G_{2k+2}(\tau)}{(2k)!}x^{2k}\end{equation} for $C=-5\G_4(\tau)$, see  Corollary \ref{diffth}.
Here $\G_{2k}(\tau)$ are  modular forms, normalized as in \cite[\S2]{Zagier}.
This allows us to write relations between the modular forms and express them
easily as polynomials in $\G_{4}(\tau)$ and $\G_{6}(\tau)$. It was noticed by Anatoly Libgober, that the equation \eqref{mojerownanie} is related to the basic equation for the Weierstrass $\wp$-function. Namely,  \eqref{weierstarss} is essentially  the regular part of the Weierstrass function
\begin{equation}\label{substitution}\wp(x)=x^{-2}+2u(x)\,,\end{equation}
see %[Stein-Shakarchi, Complex Analysis, Princeton University Press, Princeton, NJ, 2003], p. 274
see e.g.~\cite[p. 274]{SS}.
Let us differentiate
the  equation for $\wp$
$$\wp'(x)^2=4 \wp(x)^3-g_2 \wp(x)+g_3\qquad|\tfrac d{dx}\,.$$
We obtain 
%$$\wp'(x) \left(g_2+2 \wp''(x)-12 \wp(x)^2\right)=0$$
$$2 \wp''(x)\wp'(x)= \left(12 \wp(x)^2-g_2\right)\wp'(x)$$
and since  $\wp'(x)\not\equiv0$ we have
$$2 \wp''(x)=12 \wp(x)^2-g_2\,.$$
With the substitution \eqref{substitution}, dividing by 4, we arrive to \eqref{mojerownanie}, the constant agrees: 
$5\G_4(\tau)=\frac{g_2}4$.
Similar differential equation is obtained in \S\ref{flop}, equation \eqref{diffrow2}. Recursive relations for quasi-Jacobi forms follow. This is an alternative description of the  2-parameter elliptic genus, comparing with the characterization by H\"ohn \cite[Lemma 2.2.1]{Hohn}, \cite[p.~198]{HirzebruchMani}. 
\medskip

I am very grateful to Jakub Koncki for suggesting important improvements to the
first version of this paper and for many inspiring discussions. Also, I would like to thank Anatoly Libgober for enlightening    remarks  and further discussion on elliptic genus.

\section{Preliminaries and baby example}\label{baby}
Let us recall Hirzebruch construction of characteristic classes.
The {\it Hirzebruch formalism} provides a method to define {\it multiplicative
characteristic classes} of complex vector bundles using formal power series.
These classes appear naturally in various index theorems and in the definition
of genera of manifolds.
Let $E $ be a complex vector bundle over a smooth manifold $X $. A
\emph{multiplicative characteristic class} $\cc $ assigns to $E $ a cohomology
class $\cc(E) \in H^*(M; \C) $ satisfying:
$$
\cc(E \oplus F)=\cc(E) \cdot \cc(F)
$$
Such classes are defined in terms of the \emph{Chern roots} $x_1, \ldots, x_n $
of $E $.
Let $\fgf(x) \in \mathbb{Q}[[x]] $ be a formal power series with constant term
$\fgf(0)=0 $, $\fgf'(0)=1$. Given a complex vector bundle $E $ with Chern roots $x_1,
\ldots, x_n $, the w define the multiplicative
characteristic class associated to $\fgf$ by:
$$
\cc_{\fgf}(E)=\prod_{i=1}^n \frac{x_i}{\fgf(x_i)}
$$
This is a universal method to construct characteristic classes from formal
power series. For a manifold $X$ we define $\cc_\fgf(X)=\cc_\fgf(TX)$. The genus of a
compact variety $X$ associated to $\fgf$ is defined as the integral
$\int_X\cc_\fgf(X)$, or in the topological language -- the push forward to the
point.

Some classical classes and genera correspond to specific choices of $\fgf(x) $:
$$
\begin{array}{cccl}
\text{series }\fgf(x)& \text{characteritic class}& \text{genus}&\\[0.4cm]
\frac x{1+x} &\text{Chern class}& \quad \int_X c_*(TX)&\chi_{\text{top}}(X)
\\[0.4cm]
1-e^{-x} &\text{Todd class}& \quad \chi(X, \mathcal{O}_X)&\text{Todd genus}
\\[0.4cm]
{\tanh x} &\text{L-class}& \quad \sigma(X)&\text{signature} \\[0.2cm]
\text{\Large$\frac{1-e^{-(1+y)x}}{1+y e^{-(1+y)x}}$} &\text{Hirzebruch class}&
\quad \chi_y(X)&\chi_y\text{-genus}
\end{array}
$$

As a preview of the paper’s main result, let us consider characteristic
classes of smooth varieties that remain invariant under blow-ups. We make free
use of equivariant cohomology, which makes even the simplest examples nontrivial
-- and, to our surprise, decisive.
Let $bl: X'=\text{Bl}_0(\mathbb{C}^2) \to \mathbb{C}^2=X $ be the blow-up at the
origin.
The torus $\T=(\C^*)^2 $ acts on $\C^2$ coordinatwise and on the resolution.
Note that $$
H^*_\T(\C^2)=\C[x_1, x_2]\,.
$$
Here $x_1$ and $x_2$ are the weights of the coordinates of $\C^2$.
There are two fixed points in the blow-up. The corresponding weights of the
torus action are the following
$$\{ x_1,~x_2-x_1\} \qquad\text{and}\qquad \{ x_2,~x_1-x_2\}\,.$$

Our basic tool to compute equivariant push-forward is the localization theorem:
\begin{theorem}[{\cite{AB}}]\label{ABtheorem}
Let $\varphi:X\to Y$ be a $\T$-equivariant proper map of $\T$-manifolds and
$\iota_X: X^\T\to X$, $\iota_Y: Y^\T\to Y$ be the inclusions of the fixed point
sets. Denote by $e(\nu_{X^\T})$, $e(\nu_{Y^\T})$ the Euler classes of the normal
bundles.
Let $S$ be the multiplicative system generated by characters of the normal
bundles.
The following diagram commutes:
$$\begin{matrix} H^*_\T(X)&
\stackrel{e(\nu_{X^\T})^{-1}\iota_X^*}{\xrightarrow{\hspace*{1.5cm}}
}&S^{-1}H^*_\T(X^\T)\\ \\
{\varphi_*}\Bigg\downarrow\phantom{f_*}&&\phantom{f_*^\T}\Bigg\downarrow
\varphi_*^\T\\
H^*_\T(Y)&\stackrel{e(\nu_{Y^\T})^{-1}\iota_Y^*}{\xrightarrow{\hspace*{1.5cm}}
}&S^{-1}H^*_\T(Y^\T)\end{matrix}$$
If there are finitely many fixed points then the map $\varphi_*^\T$ reduces to
summation of the rational functions.
\end{theorem}

Torus localization for equvariant push forward leads to the formula:
$$\frac{bl _*\cc_\fgf(X')}{x_1x_2}=
\frac1{\fgf(x_1)\fgf(x_2-x_1)}+\frac1{\fgf(x_1-x_2)\fgf(x_2)}\in \hat
S^{-1}H^*_\T(X)=\Q[[x_1,x_2]](x_1^{-1},x_2^{-1})\,,$$
while $$\frac{\cc_\fgf(X)}{x_1x_2}= \frac1{\fgf(x_1)\fgf(x_2)}\,.
$$
The identity
$$
\frac{1}{\fgf(x_1)\fgf(x_2-x_1)}+\frac{1}{\fgf(x_2)\fgf(x_1-x_2)}-\frac{1}{\fgf(x_1)\fgf(x_2)}=0
$$
after bringing to the common denominator we obtain
$$\mathcal R=\fgf(x_2)\fgf(x_2-x_1)+\fgf(x_1)\fgf(x_2-x_1)-\fgf(x_1-x_2)\fgf(x_2-x_1)=0\in
\Q[[x_1,x_2]]$$
Let
$$\mathcal{DR}=\frac d{dx_2}\mathcal R|_{x_1=x,~x_2=0}\,.$$
Assuming $\fgf(0)=0$, $\fgf'(0)=1$ we obtain
$$\mathcal{DR}=\fgf(x) +\fgf(-x) \fgf'(x)=0\,.$$
This relation allows to inductively find the coefficients of $\fgf(x)$.
We conclude that:
\begin{corollary} If the characteristic class satisfies
$\varphi(\cc_\fgf(X_1))=\cc_\fgf(X_2)$ for a blow up $\varphi\colon X_1\to X_2$ in a
codimension two center then
$$
\fgf(x)=\frac{1-e^{-\nu x}}{\nu x} ~~\text{ for }\nu\in \Q\setminus\{0\}~~\text{ or
} ~~\fgf(x)=x
$$
i.e., $\cc_\fgf$ is proportional to the Todd class or it degenerates to a trivial
case.\end{corollary}

\begin{remark}\rm Demanding the motivic scissor relation
$$bl_*\big(\cc_\fgf(X')-\iota_*(\cc_\fgf(E))\big)=\cc_\fgf(X)-[0]\,,$$
see \cite[Corollary 0.1]{BSY},
leads to a similar equation. We find that, up to a rescaling $x$, only the
Hirzebruch class is motivic.\end{remark}

We have presented above calculation to illustrate our method in a more complex
settings, where the computations were carried out symbolically, using Wolfram
Mathematica. We provide the code in \S\ref{kody} so that readers can verify the
calculations themselves.
Further refinements include:

\medskip
-- incorporating divisor multiplicities and test for a single blow-up,

\medskip
-- restricting to special classes of singularities and test for the classical
types of flops.
\medskip

\noindent All the distinguished examples of resolutions lead to one conclusion.
Only the elliptic characteristic class with an additional factor admit a
generalization to singular
varieties.

\section{Relative characteristic classes}

Let $Y$ be a normal complex algebraic variety. It has a canonical divisor $K_Y$,
which is a Weil divisor defined up to
linear equivalence. Over the smooth locus $Y_{\text{reg}}$, this divisor
represents the line bundle of top-degree differential forms:
$$
(K_Y)|_{Y_{\text{reg}}}=\Omega^{\dim Y}_{Y_{\text{reg}}}.
$$

A variety $Y$ is said to be $\Q$-Gorenstein if a multiple of the canonical
divisor $K_Y$ is Cartier. Then it defines an element of the Picard group
$Pic(Y)\otimes \Q$, in particular it can be pulled back to a resolution of
singularities.
For a log-resolution $\psi: X \to Y$, the discrepancy divisor is defined as
$$
E=K_X-\psi^*(K_Y).
$$
If another resolution $\psi': X' \to Y$ factors through a morphism
$\varphi\colon X' \to X$, then:
$$
\varphi^*(K_X-E)=K_{X'}-E',~~~~\text{where }~~E'=K_{X'}-\psi'{}^*(K_Y).
$$
The discrepancy divisor $E$ can be uniquely written as a sum of irreducible
components of the exceptional locus of $\varphi$
$$
E=\sum_{i=1}^r a_i E_i,
$$
where $E_i$ are exceptional divisors and $a_i \in \mathbb{Q}$.

Many invariants of $Y$ are defined using combinations of Chern classes of a
resolution $c_*(TX)$, the fundamental classes $[E_i]$, and the coefficients
$a_i$. Examples include:
\begin{itemize}
\item[--] topological and Hodge-theoretic zeta functions \cite{DenefLoeser},
\item[--] elliptic genus \cite{BL}.
\end{itemize}
In a similar way the stringy Hodge numbers are defined, \cite{Batyrev}.

A study of a singular variety $Y$ boils down to working with a system of smooth
varieties $X$ with SNC (simple normal crossing) divisors. By \cite{Wlod} we can
reduce the argument to the case when the morphism between smooth varieties is a
blow-up in a smooth center intersecting cleanly the divisor.
\medskip

Let $\varphi\colon X' \to X$ be a birational morphism between smooth varieties.
Given a divisor $E \subset X$, we define its \emph{canonical pull-back}
$\varphi^{\kappa}(E)$ by the condition:
$$
K_{X'}-\varphi^{\kappa}(E)=\varphi^*(K_X-E)\,.
$$

\begin{example}\label{plasz}\rm Let $X=\C^2$ and $\varphi\colon X'\to X$ be the
blow-up at
$0$.
Let $E_i=\{z_i=0\}$ for $i=1,2$
and let $E=\alpha_1E_1+\alpha_2E_2$. Then the coefficient of the exceptional
divisor in $\varphi^\kappa(E)$ is equal to $\alpha_1+\alpha_2+1$.
\end{example}
Let $\fgf(x) \in \C[[x]]$ be a formal power series with $\fgf(0)=0$, $\fgf'(0)=1$.
Suppose $(X, E)$ is a smooth variety with a simple normal crossing $\Q$-divisor
$E=\sum_{i=1}^m a_iE_i$. We define a characteristic class of $(X,E)$ depending
on the Chern class $c_*(TX)$, the divisor classes $[E_i]$ and the coefficients
$a_i$
$$ \ccc_\fgf(X,E)\in H^*(X) \,.$$
The class $\ccc_\fgf(X, E)$ is defined as
$$
\ccc_\fgf(X,E)=\prod_{i=1}^{\dim X} x_i F_\fgf(x_i, h) \prod_{j=1}^m \frac{F_\fgf(e_j,
(a_j+1)h)}{F_\fgf(e_j,h)}\,,
$$
where
$$
F_\fgf(x,h)=\frac{\fgf(x+h)}{\fgf(x)\fgf(h)}\,.
$$
Here $x_i$ are the Chern roots and $e_i=[E_i]$.

Note that $TX=\bigoplus \mathcal{O}(E_i)$, then:
$$
\ccc_\fgf(X, E)=\prod_{i=1}^{\dim X} e_i F_\fgf(e_i, (a_i+1)h)\,.
$$
\medskip

We say that $\ccc_\fgf$ is {\it preserved by canonical transformations} if for any
birational proper morphism $\varphi\colon (X', E') \to (X, E)$ we have:
\begin{equation}\label{cantrans}
\varphi_*(\ccc(X', \varphi^{\kappa}(E)))=\ccc(X, E).
\end{equation}
To test this property it is enough to test equality for
blow-ups along smooth centers.

\begin{remark}\label{logcan} \rm Note that in order to avoid zero in the
denominator we have to assume that $a_i\neq -1$. Further taking blow-ups and
applying canonical transformations the coefficients $a_i$ can grow. Therefore to
avoid zero we assume that $a_i>-1$. Then we say that the pair $(X,E)$ has at
worst log-canonical singularities.\end{remark}

\section{Consequences of the blow-up}
We search for a conditions under which the function $\fgf$ defines a class $\ccc_\fgf$
which is invariant under canonical transformations, i.e. it satisfies
\eqref{cantrans}.

\begin{theorem} \label{pierwszy} The following are equivalent
\begin{enumerate}
\item $\ccc_\fgf$ is invariant with respect to canonical transformations.
\item Blow-up relation holds:
$$F_\fgf(x, a+b)F_\fgf(y-x, b)+F_\fgf(y, a+b)F_\fgf(x-y, a)=F_\fgf(x, a)F_\fgf(y, b)\,.$$
\item The function $\fgf(x)$ is of the form $e^{\lambda x+\mu x^2}\fgf_0(x)$ and $\fgf_0$
satisfies the diffeterential equation
$$3\fgf_0''(x)^2-4 \fgf_0'(x) \fgf_0^{(3)}(x)+\fgf_0(x) \fgf_0^{(4)}(x)+\fgf_0(x)^2
\fgf_0^{(5)}(0)=0\,,$$
$$\fgf_0(0)=\fgf_0''(0)=\fgf_0^{(3)}(0)=0\,,~~~~\fgf'(0)=1\,.$$
\end{enumerate}
Extending the coefficients to $\C$ a generic\footnote{That is an open
subset of solutions. Note that the dependence on $\tau$ is not algebraic.}
solution of the above equation is equal to:
$$\fgf(x) =e^{\lambda x+\mu x^2} \cdot\frac{\theta_\tau(\nu
x)}{\nu\,\theta'_\tau(0)}$$
where $\theta_\tau(x)$ is the Jacobi theta function, $\tau \in \mathbb{H}^+$,
$\lambda,\mu,\nu\in\C$.
\end{theorem}

\begin{remark}\rm First note that $\lambda$ in the formulation of the theorem is
completely irrelevant, since the factor
$e^{\lambda x}$ cancels in the definition of $F_\fgf$.
We will show that $\fgf(x)$ can be written as
$$\fgf(x)=x\,e^{\lambda x+\mu x^2+\mepsilon(x)}$$
with $\mepsilon(x)=\sum_{k=2}^\infty a_{2k}x^{2k}$ an even function.
Moreover the coefficients $a_{2k}$ in the expansion of $\mepsilon(x)$ for $k>3$
depend on $a_4$ and $a_6$ polynomially, $\lambda$ and $\mu$ are arbitrary.
In addition to the expression given at the last line of the theorem one has to
allow degenerate solutions, which we will discuss in \S\ref{zdegenerowane}.
\end{remark}
\medskip

{\bf Proof.}Consider the blow-up of a codimension
two center.
Locally around the center the blow up looks like a family of blow-ups of
surfaces at single points. Allowing families implies that the invariance with
respect to the
canonical transformation holds in equivariant cohomology. Then instead of an
arbitrary center of codimension 2 it is enough to consider the blowup of
$0\in\C^2$
$$\varphi\colon\ X'=\text{Bl}_0 \mathbb{C}^2 \to \mathbb{C}^2=X\,.$$
Suppose $x_1,x_2$ are the weights of $\T=(\C^*)^2$ acting on $\C^2$. Let $E_1$
and $E_2$ be the coordinate divisors: the weight of the normal bundle to $E_i$
is equal to $x_i$, as in Example \ref{plasz}. Define the divisor
$E=\alpha_1E_1+\alpha_2E_2$. The canonical pull back of $E$ has the multiplicity
$\alpha_1+\alpha_2+1$ along the exceptional divisor. The invariance of $\ccc_\fgf$
with respect to the canonical transformation is equivalent to the identity
\begin{multline*}F_\fgf(x_1, (\alpha_1+\alpha_2+2)h)F_\fgf(x_2-x_1,
(\alpha_2+1)h)+\\+F_\fgf(x_2, (\alpha_1+\alpha_2+2)h)F_\fgf(x_1 -x_2,
(\alpha_1+1)h)=\\=F_\fgf(x_1, (\alpha_1+1)h)F_\fgf(x_2, (\alpha_2+1)h)\,.
\end{multline*}
Setting the variables $x=x_1$, $y=x_2$, $a=(\alpha_1+1)h$, $b=(\alpha_2+1)h$ we
obtain the identity stated in the condition (2) of Theorem \ref{pierwszy}.
We note that the identity is satisfied for
$F_\fgf(x,h)=\frac{\fgf(x+h)}{\fgf(x)\fgf(h)}$ if and only if it is satisfied for $F_{ \fgf_0}$
with $\fgf_0(x)=e^{-\lambda x-\mu x^2}\fgf(x)$. Therefore we can substitute $\fgf$ by
$\fgf_0$ and assume that $\fgf''(0)=\fgf'''(0)=0$.
Applying the definition of $F_\fgf$ we obtain a rational expression.
We bring the terms to a common denominator and then determine the numerator:
\begin{multline*}\mathcal R_0=-\fgf\left(a+b\right) \fgf\left(x-y\right)
\fgf\left(y-x\right)
\fgf\left(a+x\right) \fgf\left(b+y\right)\\
+\fgf\left(b\right)
\fgf\left(x\right) \fgf\left(y-x\right) \fgf\left(a+x-y\right)
\fgf\left(a+b+y\right)\\+\fgf\left(a\right) \fgf\left(x-y\right)
\fgf\left(y\right) \fgf\left(a+b+x\right) \fgf\left(b-x+y\right)\,.\end{multline*}
We compute
$$
\left. \frac{d^3}{d y \,d a \,d b} \mathcal R_0\right|_{y= a=b=0}
$$
Using the initial conditions:
$$
\fgf(0)=0, \quad \fgf'(0)=1, \quad \fgf''(0)=\fgf'''(0)=0
$$
we obtain:
$$
\fgf(x)\left(\fgf(x) \fgf'(-x)+\fgf(-x) \fgf'(x)\right)=\fgf(x)\fgf(-x)^2\left(\tfrac{
\fgf(x)}{\fgf(-x)}\right)'=0\,.
$$
Hence $\fgf(x)/\fgf(-x)$ is constant, and since $\fgf(0)=0$, $\fgf'(0)=1$ we deduce that
$\fgf(x)$ is an odd function.
Let \begin{multline*}\mathcal R=\mathcal R_0/\fgf(x-y)=\fgf\left(a+b\right)
\fgf\left(x-y\right) \fgf\left(a+x\right)
\fgf\left(b+y\right)\\
~~~~~~~~~~~~~~~~ -\fgf\left(b\right) \fgf\left(x\right)
\fgf\left(a+x-y\right) \fgf\left(a+b+y\right)\\+\fgf\left(a\right)
\fgf\left(y\right) \fgf\left(a+b+x\right) \fgf\left(b-x+y\right)\,.\end{multline*}
This is a form of the Fay trisecant identity\footnote{The combinatorics of
coefficients of the Fay trisecant identity was studied by Pawe{\l} Pielasa. He
was able to show that the 7-th jet of $\fgf$ determines $\fgf$.}, \cite[Example
2.10]{RW-BS}.
Using the antisymmetry
$$
\fgf(-x)=-\fgf(x), \quad \fgf'(-x)=\fgf'(x), \quad \fgf''(-x)=-\fgf''(x), \quad \fgf'''(-x)=\fgf'''(x)
$$
and vanishing of the first four derivatives at 0, we derive a differential
equation:
let us compute
\begin{multline}\label{diffrow1}
\mathcal{DR}=\left. \frac{d^6}{d^2 y \,d^2 a \,d^2 b} \mathcal R\right|_{y= a=b=0}=\\
=2 \left( 3\fgf''(x)^2-4 \fgf'(x) \fgf^{(3)}(x)+\fgf(x) \fgf^{(4)}(x)+\fgf(x)^2 \fgf^{(5)}(0)
\right).
\end{multline}
We write the solutions $\fgf(x)$ in the form:
$$
\fgf(x)=x e^{\mepsilon(x)}\,,~~~~\mepsilon(x)=\sum_{i=2}^{\infty} a_{2i} x^{2i} \,.
$$
Denoting $d_j$ as the coefficient of $x^j$ in the series expansion of
$\mathcal{DR}$, we obtain:
$$
\begin{aligned}
& d_4=0 \\
& d_6=288 (6 a_4^2+7 a_8) \\
& d_8=720 (12 a_4 a_6+11 a_{10}) \\
& d_{10}=144 (24 a_4^3+75 a_6^2+140 a_4 a_8+143 a_{12}) \\
&\vdots
\end{aligned}
$$
The coefficient $d_j$ is a polynomial in $a_k$. It is quasi-homogeneous
of degree $j+2$, provided that we assign the weight $i$ to $a_i$. We will
show that the coefficient of $a_{j+2}$ in $d_j$ does not vanish for $j\geq 6$.
Looking at the shape of the differential equation we see that the highest
derivative is $\fgf^{(4)}$. There are also terms $\fgf'\fgf^{(3)}$, $(\fgf'')^2$ and a lower
term $\fgf^2$. Hence the coefficient of $a_{j+2}$ in $d_j$ is a polynomial of
degree at most four. To find this polynomial we list few first values for
$j\geq6$:
$$
2016, 7920, 20592, 43680, 81600, 139536, 223440, 340032, 496800, 702000, 964656, %1294560,
\dots.
$$
We check that the sequence is described by the formula $ 2(j-4)(1+j)(2+j)(3+j)$,
hence it does not vanish for $j>4$. Therefore we can
find the coefficients $a_{j}$ recursively in terms of previous ones:
{\renewcommand{\arraystretch}{1.5} $$
\begin{array}{l}
a_8 = -\frac{6}{7} a_4^2 \\
a_{10} = -\frac{12}{11} a_4 a_6 \\
a_{12} = \frac{3}{143}(32 a_4^3-25 a_6^2) \\
a_{14} = \frac{1440 }{1001}a_4^2 a_6 \\
a_{16} = -\frac{36 }{2431}(44 a_4^4-75 a_4 a_6^2) \\
a_{18} = -\frac{60 }{46189}(1392 a_4^3 a_6-275 a_6^3) \\
a_{20} = \frac{432 }{2540395}(3872 a_4^5-12125 a_4^2 a_6^2)
\end{array}
$$
}
There are no conditions for $a_4$ and $a_6$.

To sum up:
The structure of the differential equation and the expansion of the function
$\fgf(x)$ in exponential form leads to a recursive algebraic formula for the
coefficients. There are two degrees of freedom.
Additionally we can multiply $\fgf(x)$ by $e^{\mu x^2}$ with arbitrary coefficient
$\mu$. The irrelevant factor $e^{\lambda x}$ does not count.
Now the point is that we know at least one parameter family of solutions. Let
$\theta_\tau(x)$ be the Jacobi theta function.
By \cite[Theorem 3.6]{BL2} the characteristic class $\ccc_\fgf$ for
$\fgf(x)=\frac{\theta_\tau(x)}{\theta_\tau'(0)}$ is preserved by the canonical
transformation. After a brief recollection of the basic facts about modular
forms in the next section we will produce a three parameter family of solutions.
That will conclude the proof of the Theorem \ref{pierwszy}.

\section{Jacobi theta function and modular forms}
Let us recall the basic facts about Jacobi theta function. We adopt the
convention of \cite{Zagier}. Let us define the Jacobi theta function by the
triple product formula
$$\theta_\tau(x)=q^{1/8}(e^{x/2}-e^{-x/2}) \prod_{n=1}^\infty (1-q^n)(1-q^n
e^x)(1-q^n e^{-x})\,.$$
As usual $q=e^{2\pi i\tau}$, ${\rm im}(\tau)>0$.
The theta function can be written as
\begin{equation}\label{wzortheta}\theta_\tau(x)={\theta'_\tau(0)}\,x\, {\rm
Exp}\Big(-2\sum_{k=1}^\infty \frac
{\G_{k}(\tau)}{k!}x^{k}\Big)\,,\end{equation}
where $\G_{k}(\tau)$ are normalized modular forms
$$\G_k(\tau)= \begin{cases}-\frac{B_k}{2k}+\sum_{n=1}^\infty \sum_{d \mid n}
d^{k-1} q^n
&\text{ if } k\text{ is even, }\\
\phantom{-}0& \text{ if } k\text{ is odd. }\end{cases}
$$
see \cite[p.450]{Zagier}. Here $B_k$ is the Bernoulli number.

The function $\fgf(x)=\frac{\theta_\tau(x)}{\theta'_\tau(0)}$ satisfies the
condition (1) of Theorem \ref{pierwszy}.
The coefficients $a_{2k}$ appearing in the proof of Theorem \ref{pierwszy} are
of the form
$$a_{2k}=-2\frac {\G_{2k}(\tau)}{(2k)!}\,.$$
Let us quote the basic fact about  modular forms:
\begin{theorem*}[{see \cite[\S2, Prop. 4]{Zagier123}}] The modular forms
$\G_4(\tau)$
and $\G_6(\tau)$ are algebraically independent. The $\Q$-algebra generated by
$\G_4(\tau)$
and $\G_6(\tau)$ contains all $\G_k(\tau)$ for $k\geq 8$.
\end{theorem*}
Together with $\G_2(\tau)$ the forms $\G_4(\tau)$ and $\G_6(\tau)$ can be
treated as free variables.
If we replace $a_2=-\G_2(\tau)$ in $\fgf(x)$ by $-\G_2(\tau)+\mu$
the invariance of $\ccc_\fgf$ for any blow-up remains preserved since it is an
algebraic relation involving coefficients $a_{2k}$, $k\geq 1$. Together with
possibility of
rescaling $x$ we have three degrees of freedom.
The generic solution is of the form
$$\fgf(x)=e^{\mu x^2}\;\frac{\theta_\tau(\nu x)}{\nu \theta_\tau'(0)}$$
with $\mu,\nu,\tau\in \C$, $\nu\neq0$, $im(\tau)>0$.
Multiplying by the factor $e^{\lambda x}$ does not change the class $\ccc_\fgf$.
This concludes the proof of Theorem \ref{pierwszy}.
\hfill\qed

The modular form $\G_2 $, having anomalous
transformation behavior, often introduces complications in modular-invariant
constructions. However, it does not appear in the elliptic genera of Calabi--Yau
manifolds, which ensures that such genera remain weak Jacobi modular forms.
As we have shown, $\G_2 $ does not play a role in the constraints that define
admissible
functions $\fgf(x) $, further indicating its irrelevance in this context.
If necessary, $\G_2 $ can be artificially eliminated, by multiplication by the
factor $e^{\G_2(\tau)x^2}$.

\section{Differential equation}
Our proof of Theorem \ref{pierwszy} is based on the analysis of the
differential expression \eqref{diffrow1} vanishing after the substitution
$\fgf(x)=e^{\G_2 x^2}\frac{\theta_\tau(x)}{\theta_\tau'(0)}$.
Remarkably simpler differential equation is satisfied by the function $v(x)$
after the substitution
$\fgf(x)=x\,e^{-2v(x)}$:
\begin{equation}\label{rowex}12 \left(x^2 v''(x)^2+v''(x)\right)-x^2
\left(v^{(4)}(x)+5 v^{(4)}(0)\right)=0\,.\end{equation}
The function
$$v(x)=\sum_{k=2}^\infty \frac{\G_{2k}(\tau)}{(2k)!}x^{2k}\,.$$
is a solution.
We note that $v(x)$ and $v'(x)$ do not appear in the equation \eqref{rowex}.
Let us set $u(x)=v''(x)$ and assume $u(0)=0$. Finally we obtain quite elegant
form of a differential equation characterizing the sequence of the modular
forms:

\begin{corollary}\label{diffth} The function
$$u(x)=\sum_{k=1}^\infty \frac{\G_{2k+2}(\tau)}{(2k)!}x^{2k}\,.$$
satisfies the equation
$$u''(x)= 12\left(\frac{u(x)}{x^2}+u(x)^2\right)+C$$
with $C=-5\G_4(\tau)$.
\end{corollary}
As it was explained in the end of Introduction this equation can be derived directly from the basic differential equation involving the Weierstrass $\wp$-function. 
\begin{remark}\rm From the above equation we obtain a convenient induction
allowing to compute
$\G_{2k}(\tau)$ in terms of $\G_{2i}(\tau)$ for $i<k$.
A substantially different set of relations involving $\G_{2k}(\tau)$ (with a
different normalization) has was proven in \cite{MertensRolen}.

\end{remark}

\section{Bott-Samelson resolutions of Schubert varieties}
An important class of singular varieties is formed by Schubert varieties,
subvarieties of rational homogeneous spaces $G/P$. In the case of complete 
flag varieties, Schubert varieties  $X_w$ are indexed by elements of the Weyl group. 
Each has a distinguished set of resolutions, known as a Bott–Samelson resolutions, 
which is indexed by the reduced expression of the corresponding Weyl group element 
in terms of simple reflections. Any two reduced expressions of the same element 
differ by a sequence of braid relations. Therefore, verifying the independence 
of the Bott–Samelson resolution reduces to checking invariance under the braid relations. 
This condition was discussed in \cite[Equation (32-33)]{RW-BS}. The formula to be
checked is:
\begin{multline}\label{braid}F_\fgf(z_2-z_1,\mu_3-\mu_2) F_\fgf(z_3-z_2,\mu_3-\mu_1)
F_\fgf(z_2-z_1,\mu_2-\mu_1)+\\
+F_\fgf(z_1-z_2,h) F_\fgf(z_3-z_1,\mu_3-\mu_1) F_\fgf(z_2-z_1,h)=\\
=F_\fgf(z_3-z_2,\mu_2-\mu_1) F_\fgf(z_2-z_1,\mu_3-\mu_1) F_\fgf(z_3-z_2,\mu_3-\mu_2)+\\
+F_\fgf(z_2-z_3,h) F_\fgf(z_3-z_1,\mu_3-\mu_1) F_\fgf(z_3-z_2,h)\,.\end{multline}
Here $z_1,z_2,z_3,\mu_1,\mu_2,\mu_3$ and $h$ are free variables. Without a loss
of generality we can assume that $z_1=\mu_1=0$, since in the formula only the
differences $z_i-z_j$ and $\mu_i-\mu_j$ appear. The above formula is a part of
the braid relation satisfied by the elliptic Demazure-Lusztig operators,
\cite[Theorem~11.1]{W-ell}. For twisted motivic Chern classes it was discussed
in \cite[Proposition 11.9]{KW2}.
This relation already appears when checking that two resolutions of the pair
$(X_w,\partial X_w)$ give equal equivariant characteristic classes for
$w=s_1s_2s_1=s_2s_1s_2$ in the Weyl group of $SL_3(\C)$.

It is surprising that we do not obtain more solutions for $\fgf$ than in the
blow-up case. The computations are similarly elementary; however, due to the
length of the intermediate steps, we omit them here. Instead, we provide the
Wolfram Mathematica code (see \S\ref{kody}). Interestingly, the resulting
differential equation is identical to that in the blow-up case, up to
multiplication by an invertible factor.

\begin{theorem}\label{drugi}The function $\fgf$ defines a characteristic class
$\ccc_\fgf$ which does not depend on the choice of Bott-Samelson resolution of
Schubert varieties if and only if $\fgf$ satisfies the equivalent conditions of
Theorem \ref{pierwszy}.
\end{theorem}
This means that no additional characteristic classes can be defined for Schubert
varieties beyond those permitted for general Kawamata log-terminal pairs.

\section{Degenerations }\label{zdegenerowane}
Suppose $\fgf(x)=\frac{\theta_\tau(x)}{\theta'_\tau(0)}$. If $\tau \to i\infty$,
i.e., $q \to 0$.
By \cite[\S3]{Zagier} the Fourier expansion of $F_\fgf$ is following:
$$F_\fgf(x, h)=\frac{1-(ab)^{-1}}{(1-a^{-1})(1-b^{-1})}+\sum_{n=1}^\infty
\left(\sum_{d \mid n} (a^{-d} b^{-n/d}-a^d b^{n/d}) \right) q^n$$
where $a=e^x$, $b=e^h$.
As $q \to 0$, the modular form $\G_k$ converges to $\frac{B_k}{2k}$. The limit of
$F_\fgf$ we see immediately from the expansion above
$$F_\fgf(x, h) \to \frac{(1-e^{-(x+h)})}{(1-e^{-x})(1-e^{-h})}$$
This degeneration is related to the Hodge-theoretic zeta function. After a
suitable change of variables, the integral $\int_X \ccc_\fgf(X, sD)$ corresponds to
the Hodge-theoretic specialization of the global motivic zeta function
associated with
the divisor $D$. The motivic global zeta function has been defined, in
considerable generality, for $\Q$-divisors on $\Q$-Gorenstein varieties in
\cite{Veysco}.

A particularly interesting degeneration arises when we set $ q = e^{-h} \to 0 $;
see \cite{KoWe}, where the corresponding constructions are described in the
context of K-theory. Applying the Riemann--Roch transformation allows us to
transfer these classes to cohomology. In this limit, the resulting classes
depend discontinuously on the coefficients of the divisor, the dependence is
locally constant.

Another interesting degeneration cannot be obtained solely from the elliptic
characteristic class. By treating $a_2$, $a_4$, and $a_6$ as free variables, we
are free to set $a_2 = \mu$ and $a_4 = a_6 = 0$.
Then
$$\fgf(x)=x e^{\mu x^2 }, \quad x F_\fgf(x,h)=e^{2\mu h x} (1+\tfrac xh)\,.$$
Assuming that the boundary divisor is empty and $X$ is smooth the associated
characteristic classes we can be writen as
$$\ccc_\fgf(X,\emptyset)=e^{2\mu c_1(X)}\,c_*(X)$$
assuming that $h=1$. It is remarkable that this class can be extended to
singular varieties.

\section{A nice choice of generators of quasi-Jacobi forms}
In this section instead of cohomology with rational coefficients we consider
complex coefficients.
From the work of Borisov and Libgober it follows that the elliptic classes,
defined in terms of the Chern roots
$$\mathcal E\ell\ell(X)=\prod _{i=1}^{\dim X}x_i\frac
{\theta_\tau(x_i+h)\theta'_\tau(0)}{\theta_\tau(x_i)\theta_\tau(h)}$$
has a generalization to singular varieties. In particular if $\psi_i\colon
X_i\to Y$, $i=1,2$ are two crepant resolutions of a singular variety, then
$\psi_{1*}\mathcal E\ell\ell(X)=\psi_{2*}\mathcal E\ell\ell(X)$.
The elliptic genus of a variety admitting a crepant resolution belongs to the
algebra spanned by the coefficients of the expansion
$$x\, F_\tau(x,h)=x\frac
{\theta_\tau(x+h)\theta'_\tau(0)}{\theta_\tau(x)\theta_\tau(h)}=1+\sum_{i=1}^\infty
\ff_i(\tau,h)x^i\,.$$
This algebra is a polynomial algebra freely generated by $\ff_1,\ff_2,\ff_3,\ff_4$,
equal to the algebra of quasi-Jacobi forms (of depth $(k,0)$, $k>0$), see
\cite[Theorem 2.12]{Lib}. In our approach it is more natural to use another set
of generators:
\def\b#1{\widetilde E_{#1}(\tau,h)}
$$
x\, F_\tau(x,h)
=\frac x{\fgf(x)}={\rm Exp}\Big(-\sum_{i=1}^\infty \b i\Big)\,,$$
$$\begin{array}{l}
\b1=-\ff_1, \\[0.3cm]
\b2= \frac{1}{2} \left(\ff_1^2-2 \ff_2\right),\\[0.3cm]
\b3= \frac{1}{3} \left(-\ff_1^3+3 \ff_2 \ff_1-3 \ff_3\right), \\[0.3cm]
\b4= \frac{1}{4} \left(\ff_1^4-4 \ff_2 \ff_1^2+4 \ff_3 \ff_1+2 \ff_2^2-4 \ff_4\right).\\
\end{array}$$
The functions $\b k$ can be computed from the formula \eqref{wzortheta} or
directly from \cite[p. 456, (viii)] {Zagier}
$$F_\tau(x,h)= \frac {x+h}{xh}{\rm Exp}\left(-\sum_{k>0}\frac 2{k!}
\big((x+h)^k-x^k-h^k\big)\G_k(\tau)
\right)$$
(we recall that $\G_k(\tau)=0$ for $k$ odd). Thus
$$\sum_{i=1}\b i x^n= \log\Big(\frac h{h+x}\Big)+\sum_{k>0}\sum_{i=1}^{k-1}\frac
2{k!} \binom k i x^ih^{k-i}\G_k(\tau)\,,$$

\begin{equation}\label{bi}\b i=\frac {(-1)^i}i \frac1{h^i}+\sum_{k>i}\frac 2{k!}
\binom k i \G_k(\tau) h^{k-i}\,.\end{equation}
Note that $\G_2(\tau)$ appears only in $\b 1$.
As in \cite[Ch. III, eq. (10-11)]{Weil}, let us set
$$e_k(\tau)=\frac 2{(k-1)!}{\G_k(\tau)}\,,$$
$$E_i(\tau,h)=\frac1{h^i}+(-1)^i \sum_{k=1}^\infty\binom{k-1}{i-1}e_k(\tau)
h^{k-i}\,.$$
In fact, the summation is over $k\geq i$.
The equation \eqref{bi} can be transformed as follows
$$\b i=\frac {(-1)^i}i \left(\frac 1{h^i}+(-1)^i\sum_{k>i}\frac 2{(k-1)!} \binom {k-1} {i-1} \G_k(\tau)h^{k-i}\right)\,,$$
$$\b i=\frac {(-1)^i}i \Big(E_i(\tau,h)-e_i(\tau)\Big)\,.$$

\medskip

If the following section we consider a general formal power series
$$
\fgf(x) = x\, \mathrm{Exp}\Big(\sum_{i=1}^\infty b_i x^i\Big),
$$
and investigate for which sequences $ \{b_i\}_{i=1}^\infty $ the associated
characteristic class $ \cc_\fgf $ satisfies
$
\psi_{1*} \cc_\fgf(X) = \psi_{2*} \cc_\fgf(X),
$
for two crepant resolutions of a singular variety. We treat the coefficients $
b_i $ purely formally, without reference to their origin in Jacobi forms.
Ultimately, we conclude that the relations among admissible sequences $ \{b_i\}
$ coincide with those governing quasi-Jacobi forms.

As before, computations are performed in equivariant cohomology. As our initial
example of a singular space, we consider the affine cone over $ \P^1 \times \P^1
$; later, we turn to nilpotent cones.

\section{Atiyah Flop}\label{flop}

We consider two small (or in general crepant) resolutions of a given singular
variety and require that the push forward of the characteristic classes
$\cc_\fgf(X_1)$ and $\cc_\fgf(X_2)$
coincide. As a test case, we examine the simplest possible situation -- the
Atiyah flop. The following calculation was suggested in \cite{HirL}.
Let $V$ and $W$ be vector spaces of dimension 2 and let $$Y=\{A\in
\Hom(V,W)\;:\;rk(A)\leq 1\}\subset {\rm Hom}(V,W)\simeq \C^4\,.$$
There are two natural small resolutions of $Y$, determined by the general
structure of the fixed-rank locus:
$$
\begin{tikzcd}[row sep=tiny, column sep=tiny]
X_1&&Y&&X_2\\
\Hom(\mathcal O(1),W) \arrow[dd]\arrow[rr,"\psi_1"]
&&\left\{\begin{matrix}A:V\to W \\rk(
A)\leq 1\end{matrix}\right\}
&& \Hom(V,\mathcal O(-1))
\arrow[dd]\arrow[ll,,"\psi_2"']\\
&&A\arrow[dl, dashed, swap] \arrow[dr, dashed]
\\
\P(V)&\ni Ker(A) & &Im(A)\in& \P(W)
\end{tikzcd}
$$
The skew arrows are defined only on the open subset of $Y$ consisting of
matrices of the rank exactly one.

Let $\T$ be the torus of dimension 4 with basis characters $s_1,s_2,t_1,t_2$.
Suppose the weights of $\T$ acing on $V$ are $-t_1, -t_2$ and $s_1,s_2$ are the
weights of the action on $W$. There are two fixed points of the torus action at
each resolution.
The associated torus weights at the fixed points are:
\begin{align*}
&\text{Fixed point weights for }X_1 :\quad
\{t_1-t_2, s_1+t_2, s_2+t_2\}, ~~~\;
\{t_2-t_1,s_1+t_1, s_2+t_1\},
\\
&\text{Fixed point weights for }X_2 :\quad
\{s_2-s_1, s_1+t_1, s_1+t_2\}, ~~~
\{s_1-s_2, s_2+t_1, s_2+t_2\}.
\end{align*}
Localization for torus action (Theorem \ref{ABtheorem}) allows to compute the
equivariant push forward to $H^*_\T(\C^4)$:

$\bullet$ Resolution $X_1$:
$$
\mathcal S_1=\frac{ \psi_{1*}\cc_\fgf(X_1)}{e(\C^4)}=
\frac{1}{\fgf(t_1-t_2) \fgf(s_1+t_2) \fgf(s_2+t_2)}+\frac{1}{\fgf(s_1+t_1) \fgf(s_2+t_1)
\fgf(t_2-t_1)}\,.
$$

$\bullet$ Resolution $X_2$:
$$
\mathcal S_2=\frac{ \psi_{2*}\cc_\fgf(X_2)}{e(\C^4)}=
\frac{1}{\fgf(s_2- s_1) \fgf(s_1+t_1) \fgf(s_1+t_2)}+\frac{1}{\fgf(s_1-s_2) \fgf(s_2+t_1)
\fgf(s_2+t_2)}\,.
$$
We write the difference of the two localization expressions
$
\mathcal S_1-\mathcal S_2$ as a quotient, i.e.~we find the common denominator.
The numerator is equal to
\begin{multline*}
\mathcal R=\fgf(s_1-s_2) \fgf(s_2-s_1) \fgf(t_1-t_2) \fgf(s_1+t_2) \fgf(s_2+t_2) \\
+\fgf(s_1-s_2) \fgf(s_2-s_1) \fgf(t_2- t_1)\fgf(s_1+t_1) \fgf(s_2+t_1) ~~~~ \\
\phantom.\hskip80pt- \fgf(s_2- s_1) \fgf(t_1-t_2)\fgf(t_2- t_1)\fgf(s_1+t_1) \fgf(s_1+t_2) \\
-\fgf(s_1-s_2) \fgf(t_1-t_2)\fgf(t_2-t_1)\fgf(s_2+t_1) \fgf(s_2+t_2) \,.
\end{multline*}
We have assumed that $\fgf(x)$ has the form:
$$
\fgf(x)=x \exp\Big( \sum_{i=1}^{\infty} b_i x^i \Big)
$$
The relation $\mathcal R=0$ is preserved under the scaling substitution $\fgf(x)
\mapsto e^{\lambda x} \fgf(x)$, thus we can assume that $b_1=0$, hence $\fgf''(0)=0$.
Let us define the differential relation
$$\mathcal{DR}=\frac{d^4}{ds_1^2ds_2^2} \mathcal
R\Big|_{s_1=s_2=t_1=0,~t_2=t}=0\,.
$$
Having $\fgf(0)=0$, $\fgf'(0)=1$, and $\fgf''(0)=0$, the expression simplifies to:
\begin{equation}\label{diffrow2}
\mathcal{DR}= 4\fgf(-t) \left(2 \fgf'(t)^2-\fgf^{(3)}(0) \fgf(t)^2-\fgf(t)
\fgf''(t)\right)+8\fgf(t)\,.\end{equation}
Substituting the exponential form of $\fgf(x)=x\, e^{\mepsilon(x)}$, with
$\mepsilon(0)=\mepsilon'(0)=0$, we expand in powers of $t$ and extract
coefficients.
Let $d_j$ denote the coefficient of $t^j$ in the differential expression
$\mathcal{DR}$. Calculation shows that $d_j=0$ for $j<6$ and the coefficient of
$d_j$ depends on $b_{j-1}$ linearly for $j\geq 6$, for example:
\begin{align*}&d_6=-8 \left(6 b_2 b_3-5 b_5\right)\,,\\
&d_7=-\tfrac{4}{3} \left(8 b_2^3+24 b_4 b_2+27 b_3^2-42
b_6\right)\,.\end{align*}
We find the coefficients of $b_{j}$ in $d_{j+1}$:
$$
40, 56, 112, 144, 216, 264, 352, 416, 520, 600, 720, 816, 952, 1064,
1216,\dots
$$
Taking second differences of the above sequence we obtain the sequence
$$
40, -24, 40, -24, 40, -24, 40, -24, \dots
$$
This indicates that the coefficients of the original sequence split into two
quadratic subsequences. This can be easily proven from the form of the
differential equation. Alternation appears since $\fgf(-t)$ is present in
\eqref{diffrow2}, so even and odd coefficients $b_i$ should be treated
separately.
None of the coefficient vanishes, therefore one can compute $b_i$ for $i>4$.
The formula depends polynomialy on $b_2$, $b_3$ and $b_4$. Together with $b_1$
we have four degrees of freedom.
First few resulting formulas are given below:
{\renewcommand{\arraystretch}{1.5}
$$\begin{array}{l}
b_5= \frac{6}{5} b_2 b_3 \\
b_6= \frac{1}{42} \left(8 b_2^3+24 b_4 b_2+27 b_3^2\right) \\
b_7= \frac{6}{7} b_3 \left(b_2^2+b_4\right) \\
b_8= \frac{1}{21} \left(b_2^4+10 b_4 b_2^2+27 b_3^2 b_2+7
b_4^2\right) \\
b_9= \frac{1}{7} b_3 \left(4 b_2^3+12 b_4 b_2+3 b_3^2\right) \\
b_{10}= \frac{2}{1155} \left(28 b_2^5+160 b_4 b_2^3+945 b_3^2 b_2^2+340
b_4^2 b_2+540 b_3^2 b_4\right) \\
b_{11}= \frac{4}{77} b_3 \left(8 b_2^4+38 b_4 b_2^2+27 b_3^2 b_2+14
b_4^2\right) \\
\vdots\\
\end{array}$$}
We know that
$$
\fgf(x)=\frac{\theta_\tau(x+h)\, \theta'_\tau(0)}{\theta_\tau(x)\, \theta_\tau(h)}
$$
satisfies the above constraints. Modifying $\fgf(x) $ by the factor $e^{\alpha x} $
and rescaling $x $, we achieve three degrees of freedom, keeping $\tau$ fixed.

\begin{theorem}\label{AtFlop}The characteristic class $\cc_\fgf(Y)=\psi_*\cc_\fgf(X)$
does not depend on the crepant resolution $\psi:X\to Y$ if and only if
$$\fgf(x)=e^{\alpha x}\frac{\theta_\tau(\nu x+h)\theta'_\tau(0)}{\nu
\theta_\tau(\nu x)\theta_\tau(h)}$$
for $\alpha,\nu,h,\tau\in \C$, $\nu\neq0$, ${\rm im}(\tau)>0$ or $\fgf(x)$ is a
limit of functions
of the above form.
\end{theorem}

Of course the part ``if'' follows from the properties of the elliptic
characteristic class proved in \cite{BL}.

\begin{remark}\rm Comparing with \cite{Lib} we use a different expansion of
$\fgf(x)$.
The exponential expansion is convenient, because it allows easily to get rid of
the parameter $b_1$, which does not enter into the relation.\end{remark}

\section{Symplectic singularities and Grassmannian flops}
According to \cite{Beauville} a symplectic singular variety $Y$ is a normal
variety
whose regular part is a symplectic manifold and there exists a
resolution $\psi: X \to Y $ and
$\psi^*(\omega_{Y_{\text{reg}}}) $ extends over $X $ to a symplectic form.
The variety $X$ together wit the map to $Y$ is called a symplectic resolution.
The
examples of symplectic singular manifolds are closures of nilpotent orbits in a
semisimple Lie algebra, as described below.

Let us concentrate on the nilpotent orbits of the type $A$, i.e. the conjugacy
classes of nilpotent matrices. We assume that the nilpotency order is $\leq 2$.
Consider the closure of the orbit of a matrix of rank $k$ with square zero.
There are two well known symplectic resolusion, see \cite{Fu} for the
classifications of the symplectic resolutions for general case. Here both
resolutions are the cotangent
spaces of Grassmannians:
$$
\begin{tikzcd}[row sep=tiny, column sep=tiny]
X_1&&Y&&X_2\\
T^*\Grass_{n-k}(\C^n) \arrow[dd]\arrow[rr,"\psi_1"]
&&\left\{\begin{matrix}A:\C^n\to\C^n \\A^2=0,~rk A\leq
k\end{matrix}\right\}
&& T^*\Grass_{k}(\C^n)
\arrow[dd]\arrow[ll,"\psi_2"']\\
&&A\arrow[dl, dashed, swap] \arrow[dr, dashed]
\\
\Grass_{n-k}(\C^n)&\ni Ker(A) & &Im(A)\in& \Grass_{k}(\C^n)
\end{tikzcd}
$$

Let $k=1$. Both resolutions of the nilpotent cone are the isomorphic to the
cotangent
bundles of projective spaces. We consider the multiplicative class $\cc_\fgf$ and
compute the push-forward of $\cc_\fgf(T^*\mathbb{P}^{n-1})$ using
resolutions $X_1$ and $X_2$.

Let $\T=\C^*\times (\C^*)^{n}$ with coordinates $s$ and $x_i$ for $i=1,2,\dots,
n$.
The first factor acts by the scalar multiplication, and the second factor acts
by conjugation. For each $n$, define the following sums:
$$
\mathcal S_1(n) = \frac{ \psi_{1*}\cc_\fgf(X_1)}{e({\rm
Hom}(\C^n,\C^n))}=\sum_{p=1}^n \prod_{\substack{i=1 \\ i \neq p}}^n
\frac{1}{\fgf(x_i-x_p)\, \fgf(x_p-x_i+s)},
$$
$$
\mathcal S_2(n) = \frac{ \psi_{2*}\cc_\fgf(X_2)}{e({\rm
Hom}(\C^n,\C^n))}=\sum_{p=1}^n \prod_{\substack{i=1 \\ i \neq p}}^n
\frac{1}{\fgf(x_p-x_i)\, \fgf(x_i-x_p+s)}.
$$
These expressions are equal to the push forwards of the classes $\cc_\fgf(X_i)$
divided by the Euler class of $\Hom(\C^n,\C^n)$.
For example for $n=3 $:
\begin{align*}
\mathcal S_1(3)=\frac{1}{
\fgf(s+x_1-x_2)\, \fgf(-x_1+x_2)\,
\fgf(s+x_1-x_3)\, \fgf(-x_1+x_3)
} \\
+ \frac{1}{
\fgf(x_1-x_2)\, \fgf(s-x_1+x_2)\,
\fgf(s+x_2-x_3)\, \fgf(-x_2+x_3)
} \\
+ \frac{1}{
\fgf(x_1-x_3)\, \fgf(x_2-x_3)\,
\fgf(s-x_1+x_3)\, \fgf(s-x_2+x_3)
}.
\end{align*}
To obtain $\mathcal S_2(3)$ one has to change the sign of all $x$-variables. We
are looking for the functions $\fgf$ for which $\mathcal S_1(3)=\mathcal S_2(3)$.
Let
$\mathcal R(3)$ be the numerator of the difference $\mathcal S_1(3)-\mathcal
S_2(3)$ brought to a common denominator.
To study the vanishing of the difference $\mathcal R(3)$, we analyze the
coefficient of the expansion under the assumption that $\fgf(x)$ satisfies:
$$
\fgf(0)=0, \quad \fgf'(0)=1, \quad \fgf''(0)=0.
$$
We can assume that $\fgf''(0)=0$ because the condition $\mathcal R(3)$ does not
change after replacing $\fgf(x)$ by $e^{\lambda x} \fgf(x)$.
Define the differential expression:
$$
\mathcal{DR}_3 = \frac{d^9}{dx_1^4dx_2^5} \mathcal
R_3\Big|_{x_1=x_2=x_3=0}\,,
$$
{\setlength{\arraycolsep}{1pt} \begin{equation}\label{diffeqgrass}
\begin{array}{rl}
\mathcal{DR}_3=40\Big(&
360 \fgf'(s)^3 \fgf''(s) \\
-&60 \fgf(s)\fgf'(s) \left(3 \fgf''(s)^2+2 \fgf^{(3)}(0) \fgf'(s)^2+4 \fgf^{(3)}(s) \fgf'(s)\right)
\\
+&60 \fgf(s)^2\left(\fgf^{(4)}(s) \fgf'(s)+\fgf^{(3)}(s) \fgf''(s)+\fgf^{(3)}(0) \fgf'(s)
\fgf''(s)\right) \\
+&\fgf(s)^3\left(\fgf^{(5)}(0) \fgf'(s)-\fgf^{(5)}(s)\right) \\
+&\fgf(s)^4\left( \fgf^{(6)}(0)-5 \fgf^{(3)}(0) \fgf^{(4)}(0) \right)\!\Big) .\\
\end{array}
\end{equation}}
Substituting $\fgf(x)$ by $x \exp\left( \sum_{i \geq 2} b_i x^i \right)$, we
compute the Taylor coefficients of $\mathcal{DR}_3$ in powers of $s$.
The coefficient of $s^j$ are denoted $d_j$, for example:
\begin{align*}d_7&=-57600 \left( b_2^4+2 b_2^2 b_4-b_4^2-6 b_2 b_6+3 b_8 \right)
\,,\\
d_8&=-28800 \left(12 b_3 b_2^3+10 b_5 b_2^2+12 b_3 b_4 b_2-42 b_7 b_2-10 b_4
b_5-18 b_3 b_6+27 b_9\right)
\end{align*}
and so on. The coefficient $d_j$ is linear in the variable $b_{j+1}$, and $d_j$
does not depend on $b_k$ for $k>j+1$.
We compute the sequence of the coefficient of $b_{j+1}$ in $d_j$ for $j\geq 7$.
Since the number differentiations in the monomials of \eqref{diffeqgrass} is at
most 5 (there are monomials $(\fgf')\fgf''$,~ $\fgf'(\fgf'')^2$, \dots) the coefficient of
$b_{j+1}$ in $d_j$ is a polynomial of degree $\leq 5$.
The first few coefficients are equal:
$$
-172800, -777600, -2217600, -5068800, -10108800, -18345600,
-31046400, -49766400,\dots$$
Having initial values we find that
the coefficient of $b_{j+1}$ is equal to
$$-240 (j-6) (j-5) (j-2) (j+1) (j+2),$$ hence it does not vanish for $j>6$.
Therefore the unknown variables $b_j$ for $j>7$ is determined by $b_i$ for
$i=2,3,\dots,7$.
\bigskip

A similar analysis for $n=4$ yields a higher-order differential equation. Define
$$
\mathcal{DR}_4 = \frac{d^{15}}{dx_1^5dx_2^5dx_3^5}\mathcal R_4\Big|_{x_1=
x_2=x_3= x_4=0}.
$$
This leads to a complicated expression involving derivatives of $\fgf(x)$ up to
sixth order, and again leads to relations among the coefficients $b_i$ in the
expansion of $\fgf(x)$. The first nonvanishing coefficients of $d_j$ is
$$d_{13}=-27648000 \left(128 b_2^5+160 b_4 b_2^3-780 b_6 b_2^2-240 b_4^2 b_2+600
b_8 b_2+240 b_4 b_6-165 b_{10}\right)\,.$$
Having only the relations for $n=4$ we would be able only to determine the
variables starting from $b_{10}$ as polynomials in $b_i$, $i=2,3,\dots,9$.
We claim that together with the relations deduced from $n=3$ all the
coefficients $b_j$ can be expressed as polynomials in $b_2$, $b_3$, $b_4$. This
is checked by a direct computation in Wolfram Mathematica. See \S\ref{kody}.
\bigskip

We have analyzed two different Grassmannian flops in type A via torus
localization and found two differential equation. The differential equations
produce polynomial constraints among the coefficients of the generating power
series $\fgf(x)$.
Having two independent differential equations involving the function $\fgf$ we find
that there is only four degrees of freedom.

\begin{theorem}\label{czwarty}Suppose that $\cc_\fgf$ is invariant with respect to
Grassmannian flops, then $\fgf$ is of the form as in Theorem \ref{AtFlop}.
\end{theorem}

\begin{remark}\rm An examination of the nilpotent orbit of $5 \times 5$ matrices
of rank at most 2 and square zero leads to the same constraints on the
coefficients as in the case of $3 \times 3$ matrices of rank one. There are many
other tractable cases that one could compute, but we have not studied them
systematically, since just two simple Grassmannian flops already determine the
elliptic genus.
\end{remark}

\section{Codes}\label{kody}
This script is written in Wolfram Mathematica. It contains all the calculations
discussed in the paper. The source:
\href{https://www.mimuw.edu.pl/~aweber/4paths}{\color{blue}\tt
https://www.mimuw.edu.pl/\~{}aweber/4paths}

%\begin{lstlisting}[language=Mathematica]

\begin{lstlisting}[language=Mathematica]
md=20; (* maximal defree *)
expansion=Normal[Series[u Exp[Sum[a[i]u^i,{i,4,md+5,2}]],{u,0,md+5}]];
ff[x_]:=expansion/.u->x
FF[x_,y_]:=f[x+y]/(f[x]f[y])
Print["Blow-up relation"]
LHS=F[x1,b1+b2] F[x2-x1,b2]+F[x2,b1+b2] F[x1-x2,b1];
RHS=F[x1,b1] F[x2,b2];
Rel0=Numerator[Factor[LHS-RHS/.F->FF]]
Print["if f''[0]=0, then f is an odd function"]
inic={f[0]->0,f'[0]->1,f''[0]->0};
Factor[D[Rel0,{x2,1},{b1,1},{b2,1}]/.{x2->0,b1->0,b2->0,x1->x}/.inic]==0
Print["relation, assuminf that f is antisymmetric"]
Rel=Factor[Rel0/f[x1-x2]/.f[x2-x1]->-f[x1-x2]]
Print["differential relation for odd function"]
inic={f[0]->0,f'[0]->1,f''[0]->0,f'''[0]->0};
asym={f[-x]->-f[x],f'[-x]->f'[x],f''[-x]->-f''[x],f'''[-x]->f'''[x],f''''[-x]->-f''''[x]};
DR=Factor[D[Rel,{x2,2},{b1,2},{b2,2}]/.{x2->0,b1->0,b2->0,x1->x}/.inic/.asym]
Print["the sequence of relations, d[j] depends linearly on a[j+2]"]
Do[d[j]=Factor[Coefficient[DR/.f->ff,x,j]],{j,0,md,2}]
Do[Print["d[",j,"]=",d[j]],{j,0,16,2}]
Print["coefficients obey polynomial pattern -2(j-4)(1+j)(2+j)(3+j)"]
Table[-2(j-4)(1+j)(2+j)(3+j)==Coefficient[d[j],a[j+2]],{j,6,md-2,2}]
Print["expressions for a[j], j>8"]
p={};Do[p=Union[p,Solve[0==d[i]/.p,a[i+2]][[1]]],{i,6,md,2}];
Column[Expand[p]]

Print["a simpler form of the equation"]
inic={v[0]->0,v'[0]->0,v''[0]->0};
Expand[Solve[0==DR/.f->Function[x,x  E^(-2 v[x])]/.inic,v''''[x]]]

Print["Braid relation"]
LHS=F[z2-z1,m3-m2] F[z3-z2,m3-m1] F[z2-z1,m2-m1]+F[z1-z2,h] F[z3-z1,m3-m1] F[z2-z1,h];
RHS=F[z3-z2,m2-m1] F[z2-z1,m3-m1] F[z3-z2,m3-m2]+F[z2-z3,h] F[z3-z1,m3-m1] F[z3-z2,h];
Rel0=Numerator[Factor[LHS-RHS/.F->FF/.m1->0/.z1->0]]
Print["if f''[0]=0,then f is an odd function"]
inic={f[0]->0,f'[0]->1,f''[0]->0};
Factor[D[Rel0,{m2,1},{m3,1},{z2,3},{h,1}]/.{m2->0,m3->0,z2->0,h->0}/.inic]==0
Print["relation,assuminf that f is antisymmetric"]
Rel=Factor[Rel0/(f[z2]f[z2-z3])/.{f[-z2]->-f[z2],f[z3-z2]->-f[z2-z3]}]
Print["differential relation for odd function"]
inic={f[0]->0,f'[0]->1,f''[0]->0,f'''[0]->0};
asym={f[-x]->-f[x],f'[-x]->f'[x],f''[-x]->-f''[x],f'''[-x]->f'''[x],f''''[-x]->-f''''[x]};
DBr=Factor[D[Rel,{m2,1},{m3,2},{z2,2},{h,4}]/.{m2->0,m3->0,z2->0,z3->x,h->0}/.inic/.asym]
Print["new differential equation differs by an invertible factor"]
Factor[DBr/DR]

Print["Atiyah Flop"]
expansion=Normal[Series[u Exp[Sum[b[i]u^i,{i,2,md+14}]],{u,0,md+14}]];
ff[x_]:=expansion/.u->x
Print["the relation"]
SA1=1/(f[t1-t2] f[s1+t2] f[s2+t2])+1/(f[s1+t1] f[s2+t1] f[-t1+t2]);
SA2=1/(f[-s1+s2] f[s1+t1] f[s1+t2])+1/(f[s1-s2] f[s2+t1] f[s2+t2]);
RA=Numerator[Factor[SA1-SA2]]
Print["differential equation assuminf f(0)=0, f'(0)=1, f''(0)=0"]
inic={f[0]->0,f'[0]->1,f''[0]->0};
DR=Factor[D[RA,{s1,2},{s2,2}]/.{s1->0,s2->0,t1->0,t2->t}/.inic];
Collect[DR,f[-t],Factor]
Print["the sequence of relations: d[j] depends linearly on b[j-1]"]
Do[d[j]=Factor[Coefficient[DR/.f->ff,t,j]],{j,0,md}]
Do[Print["d[",j,"]=",d[j]],{j,0,10}]
Print["coefficient obey a polynomial pattern, 4(j-4)(j-1), j even "]
Table[Coefficient[d[j],b[j-1]]==4(j-4)(j-1),{j,6,md,2}]
Print["coefficient obey a polynomial pattern, 4(j-4)j, j odd "]
Table[Coefficient[d[j],b[j-1]]==4(j-5)j,{j,5,md,2}]
Print["expressions for a[j] for j>4"]
p={};Do[p=Union[p,Solve[0==d[i]/.p,b[i-1]][[1]]],{i,6,md}];
Column[Expand[p]]

Print["frassmannian flop"]
S1[n_]:=Sum[Product[1/(f[x[i]-x[p]]f[x[p]-x[i]+s]),{i,Complement[Range[n],{p}]}],{p,1,n}]
S2[n_]:=Sum[Product[1/(f[x[p]-x[i]]f[x[i]-x[p]+s]),{i,Complement[Range[n],{p}]}],{p,1,n}]
R[n_]:=Numerator[Factor[S1[n]-S2[n]]]
inic={f[0]->0,f'[0]->1,f''[0]->0};

Print["differential equation for the nilpotent cone, n=3"]
DR3=Factor[D[Factor[D[R[3]/.x[1]->0,{x[2],5}]/.x[2]->0]/.inic,{x[3],4}]/.x[3]->0/.inic];
Collect[DR3/40,f[s],Factor]
Print["the sequence of relations, d[j] depends linearly on b[j-1]"]
Do[d[j]=Factor[SeriesCoefficient[DR3/.f->ff,{s,0,j}]],{j,0,md}]
Do[Print["d[",j,"]=",d[j]],{j,0,16}]
Print["the polynomial pattern -240(j-6)(j-5)(j-2)(j+1)(j+2)"]
Table[Coefficient[d[j],b[j+1]]==-240(j-6)(j-5)(j-2)(j+1)(j+2),{j,7,md}]
Print["expressions for b[j] for j>7"]
p3={};Do[p3=Union[p3,Solve[0==d[i]/.p3,b[i+1]][[1]]],{i,7,md}];
Column[Expand[p3]]

Print["differential equation for the nilpotent cone, n=4"]
DR4=Factor[D[R[4]/.x[1]->0,{x[2],5}]/.x[2]->0/.inic];
DR4=Factor[D[DR4,{x[3],5}]/.x[3]->0/.inic];
DR4=Collect[D[DR4,{x[4],5}]/.x[4]->0/.inic,f[s],Factor];
Collect[DR4/144000,f[s]]
Print["the sequence of relations, d[j] depends linearly on b[j-3]"]
Do[d[j]=Factor[SeriesCoefficient[DR4/.f->ff,{s,0,j}]],{j,0,md}]
Do[Print["d[",j,"]=",d[j]],{j,0,16}]
Print["polynomial pattern 144000(j-12)(j-11)(j-10)(j-3)(j-2)(7j-43)"]
Table[Coefficient[d[j],b[j-3]]==144000(j-12)(j-11)(j-10)(j-3)(j-2)(7j-43),{j,13,md}]
Print["expressions for b[j] j>9"]
p4={};Do[p4=Union[p4,Solve[0==d[i]/.p4,b[i-3]][[1]]],{i,13,md}];
Column[Expand[p4]]

Print["combined relations"]
diff=Factor[(Table[b[i],{i,10,13}]/.p3)-Table[b[i],{i,10,13}]/.p4/.p3];
Solve[diff=={0,0,0,0},Table[b[i],{i,5,7}]]
\end{lstlisting}

%\end{lstlisting}

The code was intentionally written to be as simple as possible. We did not
implement optimizations for computational efficiency. Still the code performs
well for a reasonable range of degrees.

%\bibliographystyle{plain}
%\bibliography{nothing}

\end{document}